\documentclass{amsart}
\usepackage{amssymb}

\usepackage{amsmath}
\usepackage{amscd}
\usepackage{thmdefs}

\input{tcilatex}
\theoremstyle{definition}
\theoremstyle{remark}
\numberwithin{equation}{section}

\input tcilatex

\begin{document}
\title[Operator-Valued Moments of the Generating Operator of $L(F_{2})$]{Operator-Valued Moment Series of the Generating Operator of $L(F_{2})$ over
the Commutator Group von Neumann Algebra $L(K)$}
\author{Ilwoo Cho}
\address{Dep. of Math, Univ. of Iowa, Iowa City, IA, U. S. A}
\email{ilcho@math.uiowa.edu}
\maketitle

\begin{abstract}
In this paper, we will consider the generating operator $x=a+b+a^{-1}+b^{-1}$
of the free group factor $L(F_{2}),$ where $F_{2}=\,<a,b>$ is the free group
with two generators $a$ and $b.$ Let $K=\,<h>$ be the commutator group of $a$
and $b$ with one generator $h=aba^{-1}b^{-1}.$ Then we can construct the
group von Neumann algebra $L(K)$ and the conditional expectation $%
E:L(F_{2})\rightarrow L(K),$ defined by $E\left( \underset{g\in F_{2}}{\sum }%
\alpha _{g}g\right) =\underset{k\in K}{\sum }\alpha _{k}k,$ for all $%
\underset{g\in F_{2}}{\sum }\alpha _{g}g$ in $L(F_{2}).$ Then $\left(
L(F_{2}),E\right) $ is the $W^{*}$-probability space with amalgamation over $%
L(K).$ In this paper, we will compute the \textbf{trivial} operator-valued
moment series of the generating operator $a+b+a^{-1}+b^{-1}$ of $L(F_{2}),$
over $L(K).$ This computation is the good example for studying the
operator-valued distribution, since the operator-valued moment series of
operator-valued random variables contain algebraic and combinatorial free
probability information about the opeartor-valued distributions.
\end{abstract}

\strut

From mid 1980's, Free Probability Theory has been developed. Here, the
classical concept of Independence in Probability theory is replaced by a
noncommutative analogue called Freeness (See [9]). There are two approaches
to study Free Probability Theory. One of them is the original analytic
approach of Voiculescu and the other one is the combinatorial approach of
Speicher and Nica (See [1], [2] and [3]). \medskip Speicher defined the free
cumulants which are the main objects in Combinatorial approach of Free
Probability Theory. The free cumulants of random variables are gotten from
the free moments of random variables via the M\"{o}bius inversion (or vice
versa). But in this paper, we will concentrate only on computing the free
moments of certain random variables (See [3]). Also, Speicher considered the
operator-valued free probability theory, which is also defined and observed
originally by Voiculescu (See [9]). In this paper, we will observe the
important example of such operator-valued free probability.

\strut

Let $F_{N}$ be a free group with $N$-generators and let $L(F_{N})$ be the
free group factor defined by

\strut

\begin{center}
$L(F_{N})=\overline{\mathbb{C}[F_{N}]}^{w}.$
\end{center}

\strut

In this paper, by using so-called the recurrence diagram found in [13] and
[14], we will compute the \textbf{trivial} operator-valued moment series of
the genarating operator $G$ of the free group factor $L(F_{2}),$ defined by

\strut

\begin{center}
$G=a+b+a^{-1}+b^{-1}\in L(F_{2}),$
\end{center}

\strut 

over the group von Neumann algebra $L(K),$ where $F_{2}=\,<a,b>$ and $%
K=\,<aba^{-1}b^{-1}>$ is the commutator group of $a$ and $b$ in $F_{2}.$ i.e,

\strut 

\begin{center}
$L(K)\overset{def}{=}\overline{\Bbb{C}[<aba^{-1}b^{-1}>]}^{w}.$ 
\end{center}

\strut 

Throughout this paper, we will fix $a$ and $b$ as the generators of the free
group $F_{2}$ and we will also fix $h=aba^{-1}b^{-1}$ as the generator of
the group $K$ which is group isomorphic to the integers $\mathbb{Z}.$ Notice
that the commutator group $K$ is a subgroup of $F_{2}$ and hence the group
von Neumann algebra $L(K)$ is a $W^{*}$-subalgebra of the free group factor $%
L(F_{2}).$ Let $x$ be an operator in $L(F_{2}).$ Then there exists the
Fourier expansion of $x,$

\strut

\begin{center}
$x=\underset{g\in F_{2}}{\sum }\alpha _{g}u_{g},$ \ with $\alpha _{g}\in %
\mathbb{C},$ for all $g\in F_{2}.$
\end{center}

\strut

We can regard all $g\in F_{2}$ as unitaries $u_{g}$ in $L(F_{2}).$ For the
convenience, we will denote these unitaries $u_{g}$ just by $g.$ With this
notation, it is easy to check that

$\strut $

\begin{center}
$g^{*}=u_{g}^{*}=u_{g}^{-1}=u_{g^{-1}}=g^{-1}$ in $L(F_{2}),$
\end{center}

\strut

where $g^{-1}$ is the group inverse of $g$ in $F_{2}.$ We can define the
conditional expectation $E:L(F_{2})\rightarrow L(K)$ by

\strut

\begin{center}
$E\left( \underset{g\in F_{2}}{\sum }\alpha _{g}g\right) =\underset{k\in K}{%
\sum }\alpha _{k}k.$
\end{center}

\strut

Then we have the $W^{*}$-probability space $\left( L(F_{2}),E\right) $ with
amalgamation over $L(K).$ Let $G$ be the generating operator, which is also
said to be the radial operator of $F_{2}$, $a+b+a^{-1}+b^{-1}$ in $L(F_{2}).$

\strut 

Now, let $B\subset A$ be von Neumann algebras and assume that there is a
conditional expectation $F:A\rightarrow B.$ Then the algebraic pair $\left(
A,F\right) $ is called the $W^{*}$-probability space with amalgamation over $%
B$ (in short, the $W^{*}$-probability space over $B$). We say that each
element $a$ of $\left( A,F\right) $ is a $B$-valued (or operator-valued)
random variable. Voiculescu defined the $n$-th $B$-valued (or
operator-valued) moment of the $B$-valued random variable $a$ by

\strut 

\begin{center}
$F\left( b_{1}a...b_{n}a\right) ,$ for all $n\in \Bbb{N},$
\end{center}

\strut 

where $b_{1},...,b_{n}\in B$ are arbitrary. (Recall that, in [1], Speicher
defined the $n$-th moment of $a\in (A,F)$ by $F\left(
ab_{2}a...b_{n}a\right) ,$ for arbitrary $b_{2},...,b_{n}\in B$) Different
from the scalar-valued case when $B=\Bbb{C},$ since $B$ and $A$ are not
cummutative, in general, we have to consider the case when $b_{1},...,b_{n}$
are not $1_{B}=1_{A}.$ When $b_{1}=...=b_{n}=1_{B},$ we say that the $n$-th $%
B$-valued moment of $a$ is the trivial $n$-th $B$-valued moment of $a.$
However, if $B$ and $A$ are commutative (in other words, if $B=C_{A}(B)$),
then the trivial moment $F(a^{n})$ contains the full free probability data
for the $B$-valued moments, since

\strut 

\begin{center}
$F(b_{1}a...b_{n}a)=\left( b_{1}\cdot \cdot \cdot b_{n}\right) F(a^{n}),$
for all $n\in \Bbb{N}.$
\end{center}

\strut 

Remark that, in our case, $L(K)$ and $L(F_{2})$ are not commutative,
because, for instance,

\strut 

\begin{center}
$ah=a^{2}ba^{-1}b^{-1}\neq aba^{-1}b^{-1}a=ha$.
\end{center}

\strut 

In this paper, we will concentrate on computing the trivial $L(K)$-valued
moments of the generating operator $G=a+b+a^{-1}+b^{-1}$ of $L(F_{2}).$ As
we mentioned in the previous paragraph, since $L(K)$ and $L(F_{2})$ are not
commutative, these trivial operator-valued moments of $G$ do not contain the
complete free probability data of $G.$ But the computation, itself, provides
the way how to compute the operator-valued moments of $G$.

\strut \strut \strut 

It is easy to see that the first, second and third trivial $L(K)$-moments of 
$G$ vanish. i.e,

\strut

\begin{center}
$E\left( G^{k}\right) =0_{L(K)},$ for $k=1,2,3,$
\end{center}

\strut

since $G^{k}$ does not contain the $h^{n}$-term, for $k=1,2,3$ and for $n\in %
\mathbb{Z},$ where $h=aba^{-1}b^{-1}$ and $h^{-1}=bab^{-1}a^{-1}.$ However,
fourth trivial $L(K)$-moment $E(G^{4})$ of $G$ contains the $h$-term and the 
$h^{-1}$-term. So, finding the trivial $L(K)$-moments of $G$ is to find the $%
h^{k}$-terms of $G^{n},$ for all $k\in \mathbb{Z}$ and $n\in \mathbb{N}.$

\strut

The following recurrence diagram will play a key role to find such trivial
operator-valued moment series of the generating operator $G$ of $L(F_{2})$ ;

\strut

\strut

\begin{center}
$
\begin{array}{llllllllllll}
&  &  &  &  &  &  &  &  &  & p_{0}^{2} & =2N \\ 
&  &  &  &  &  &  &  &  &  & \downarrow &  \\ 
&  &  &  &  &  &  &  &  &  & q_{1}^{3} & =(2N-1)+2N \\ 
&  &  &  &  &  &  &  &  & \swarrow \swarrow & \searrow \searrow &  \\ 
&  &  &  &  &  &  &  & p_{2}^{4} &  &  & p_{0}^{4} \\ 
&  &  &  &  &  &  & \swarrow \swarrow & \searrow &  & \swarrow &  \\ 
&  &  &  &  &  & q_{3}^{5} &  &  & q_{1}^{5} &  &  \\ 
&  &  &  &  & \swarrow \swarrow & \searrow &  & \swarrow &  & \searrow
\searrow &  \\ 
&  &  &  & p_{4}^{6} &  &  & p_{2}^{6} &  &  &  & p_{0}^{6} \\ 
&  &  & \swarrow \swarrow &  & \searrow &  & \swarrow & \searrow &  & 
\swarrow &  \\ 
&  & q_{5}^{7} &  &  &  & q_{3}^{7} &  &  & q_{1}^{7} &  &  \\ 
& \swarrow \swarrow &  & \searrow &  & \swarrow &  & \searrow & \swarrow & 
& \searrow \searrow &  \\ 
p_{6}^{8} &  &  &  & p_{4}^{8} &  &  & \text{ \ \ \ }p_{2}^{8} &  &  &  & 
p_{0}^{8} \\ 
\vdots &  &  &  & \vdots &  &  & \text{ \ \ \ }\vdots &  &  &  & \vdots
\end{array}
$
\end{center}

\strut

where

\begin{center}
$\swarrow \swarrow $ \ : \ $(2N-1)+[$former term$]$

$\searrow $ \ \ \ \ : \ \ $(2N-1)\cdot [$former term$]$

$\swarrow $ \ \ \ \ : \ \ $\cdot +[$former term$]$
\end{center}

and

\begin{center}
$\searrow \searrow $ \ : \ \ $(2N)\cdot [$former term$].$
\end{center}

\strut

In this paper, we obtain good applications about the above recurrence
diagram. We will re-compute the moment series of the generating operators of
the free group factor $L(F_{N}),$ for all $N\in \mathbb{N},$ by using the
above recurrence diagram. This would be the one application of this
recurrence diagram (See Chapter 1). When $N=2,$ we can apply this recurrence
diagram to compute the trivial operator-valued moment series of the
generating operator of $L(F_{2})$ (See Chapter 2). Remark that to study
(scalar-valued or operator-valued) moment series of elements in an operator
algebra is to study (scalar-valued or operator-valued) free distributions of
elements in that operator algebra. So, the computations in this paper about
generating operators contain the free probability information about free
distribution of those generating operators. And the free probability
information is determined by the above recurrence diagram.

\strut

In Chapter 1, we will re-compute the (scalar-valued) moment series of the
generating operator of $L(F_{N}),$ by using the recurence diagram found in
[13] and [14]. The moment series of the generating operator of $L(F_{N})$ is
already known, but here we will compute it again, by using the above
recurrence diagram. In Chapter 2, by using the reccurence diagram when $N=2$%
, we will compute the trivial operator-valued moment series of the
generating operator $G=a+b+a^{-1}+b^{-1}$ of $L(F_{2}).$ Remark that the
moment series in Chapter 1 is a scalar-valued moment series and the trivial
operator-valued moment series in Chapter 2 is a operator-valued
(operator-valued) moment series.

\strut

\strut

\strut

\section{Moment Series of the Generating Operator of $L(F_{N})$}

\strut

\strut

Let $A$ be a von Neumann algebra and let $\tau :A\rightarrow \mathbb{C}$ be
the normalized faithful trace. Then we call the algebraic pair $(A,\tau ),$
the $W^{*}$-probability space and we call elements in $(A,\tau ),$ random
variables. Define the collection $\Theta _{s},$ consists of all formal
series without the constant terms in noncommutative indeterminants $%
z_{1},...,z_{s}$ ($s\in \mathbb{N}$). Then we can regard the moment series
of random variables as elements of $\Theta _{s}.$

\strut

\begin{definition}
Let $(A,\tau )$ be a $W^{*}$-probability space with its normalized faithful
trace $\tau $ and let $a\in \left( A,\tau \right) $ be a random variale. The
moment series of $a$ is defined by the formal series in $\Theta _{1},$

\strut 

$\ \ \ \ \ \ \ \ \ M_{a}(z)=\sum_{n=1}^{\infty }\tau (a^{n})\,z^{n}$ .

\strut 

The coefficients $\tau (a^{n})$ are called the $n$-th moments of $a,$ for
all $n\in \mathbb{N}.$
\end{definition}

\strut \strut

Let $H$ be a group and let $L(H)$ be a group von Neumann algebra. i.e,

\strut

\begin{center}
$L(H)=\overline{\mathbb{C}[H]}^{w}.$
\end{center}

\strut

Precisely, we can regard $L(H)$ as a weak-closure of the group algebra
generated by $H$ and hence

\strut

\begin{center}
$L(H)=\overline{\{\underset{g\in H}{\sum }t_{g}g:g\in H\}}^{w}.$
\end{center}

\strut

It is well known that $L(H)$ is a factor if and only if the given group $H$
is icc. (Since the free group $F_{N}$ with $N$-generators is icc, the von
Neumann group algebra $L(F_{N})$ is a factor and it is called the free group
factor.)

\strut

Now, define the canonical trace $\tau :L(H)\rightarrow \mathbb{C}$ by

\strut

\begin{center}
$\tau \left( \underset{g\in H}{\sum }t_{g}g\right) =t_{e_{H}},$ \ for all \ $%
\underset{g\in H}{\sum }t_{g}g\in L(H),$
\end{center}

\strut

where $e_{H}$ is the identity of the group $H.$ It is easy to check that the
trace $\tau $ is normalized and faithful. So, the algebraic pair $\left(
L(H),\tau \right) $ is a $W^{*}$-pobability space$.$ \strut Assume that the
group $H$ has its generators $\{g_{j}\,:\,j\in I\}.$ We say that the operator

\strut

\begin{center}
$G=\underset{j\in I}{\sum }g_{j}+\underset{j\in I}{\sum }g_{j}^{-1},$
\end{center}

\strut

the generating operator of $L(H).$ For instance, if we have a free group $%
F_{N}=\,<g_{1},...,g_{N}>,$ then the generating operator of the free group
factor $L(F_{N})$ is

\strut

\begin{center}
$g_{1}+...+g_{N}+g_{1}^{-1}+...+g_{N}^{-1}.$
\end{center}

\strut

Rest of this chapter, we will consider the moment series of the generating
operator $G$ of $L(F_{N}).$

\strut

From now, fix $n\in \mathbb{N}.$ And we will denote free group factor $%
L(F_{N})$ by $A$. i.e

\strut

\begin{center}
$A=\overline{\{\underset{g\in F_{N}}{\sum }t_{g}g:t_{g}\in \mathbb{C}\}}%
^{w}. $
\end{center}

\strut

Recall that there is the canonical trace $\tau :A\rightarrow \mathbb{C}$
defined by

\strut

\begin{center}
$\tau \left( \underset{g\in F_{N}}{\sum }t_{g}g\right) =t_{e},$
\end{center}

\strut

where $e\in F_{N}$ is the identity of $F_{N}$ and hence $e\in L(F_{N})$ is
the unity $1_{L(F_{N})}.$ The algebraic pair $\left( L(F_{N}),\tau \right) $
is a $W^{*}$-probability space. Let $G$ be the generating operator of $%
L(F_{N}).$ i.e

\strut

\begin{center}
$G=g_{1}+...+g_{N}+g_{1}^{-1}+...+g_{N}^{-1},$
\end{center}

\strut

where \ $F_{N}=\,<g_{1},...,g_{N}>.$ It is well-known that if we denote the
sum of all words with length $n$ in $\{g_{1}$ $,$ $g_{1}^{-1}$ $,$ $...,$ $%
g_{N},$ $g_{N}^{-1}\}$ by

$\strut $

\begin{center}
$X_{n}=\underset{\left| w\right| =n}{\sum }w\in A$, for all $n\in \mathbb{N}%
, $
\end{center}

\strut

then the following two recurrence relations (1.1) and (1.2) hold true,
whenever $N\geq 2$ ;

\strut

(1.1) $\ \ \ \ \ \ \ \ \ X_{1}X_{1}=X_{2}+2N\cdot e$ \ \ \ \ ($n=1$)

\strut

\ \ \ \ \ \ and

\strut

(1.2) $\ \ \ \ \ \ X_{1}X_{n}=X_{n+1}+(2N-1)X_{n-1}$ $\ \ \ \ (n\geq 2)$

\strut

(See [15]). In our case, we can regard our generating operator $G$ as $X_{1}$
in $A,$ by the very definition of $G.$

\strut

By using the relation (1.1) and (1.2), we can express $G^{n}$ in terms of $%
X_{k}$'s ; For example, $G=X_{1},$

\strut

$G^{2}=X_{1}X_{1}=X_{2}+2N\cdot e,$

\strut

$G^{3}=X_{1}\cdot X_{1}^{2}=X_{1}\left( X_{2}+(2N)e\right)
=X_{1}X_{2}+(2N)X_{1}$

$\ \ \ \ =X_{3}+(2N-1)X_{1}+(2N)X_{1}=X_{3}+\left( (2N-1)+2N\right) X_{1},$

\strut

$G^{4}=X_{4}+\left( (2N-1)+(2N-1)+2N\right) X_{2}+(2N)\left(
(2N-1)+(2N)\right) e,$

\strut

$G^{5}=X_{5}+\left( (2N-1)+(2N-1)+(2N-1)+2N\right) X_{3}$

\strut

$\ \ \ \ \ \ \ \ \ +\left( (2N-1)\left( (2N-1)+(2N-1)+(2N)\right)
+(2N)\left( (2N-1)+(2N)\right) \right) X_{1},$

\strut

$G^{6}=X_{6}+\left( (2N-1)+(2N-1)+(2N-1)+(2N-1)+2N\right) X_{4}$

\strut

$\,\,\,\ \ \ \ \ \ \ \ \ +\{(2N-1)\left( (2N-1)+(2N-1)+(2N-1)+(2N)\right) $

$\ \ \ \ \ \ \ \ \ \ \ \ \ \ \ \ \ \ \ \ \ \ \ \ \ \ \ \ \ \ \ \
+(2N-1)\left( (2N-1)+(2N-1)+(2N)\right) $

$\ \ \ \ \ \ \ \ \ \ \ \ \ \ \ \ \ \ \ \ \ \ \ \ \ \ \ \ \ \ \ \ \ \ \ \ \ \
\ \ \ \ \ \ \ \ \ \ \ \ \ \ +(2N-1)((2N-1)+(2N))\}X_{2}$

\strut

$\ \ \ \ \ \ \ \ \ \ +(2N)\left( (2N-1)\left( (2N-1)+(2N-1)+(2N))\right)
+(2N)((2N-1)+(2N))\right) e,$

etc.

\strut

So, we can find a recurrence relation to get $G^{n}$ ($n\in \mathbb{N}$)
with respect to $X_{k}$'s ($k\leq n$). Inductively, $G^{2k-1}$ and $G^{2k}$
have their representations in terms of $X_{j}$'s as follows ;

\strut

\begin{center}
$%
G^{2k-1}=X_{1}^{2k-1}=X_{2k-1}+q_{2k-3}^{2k-1}X_{2k-3}+q_{2k-5}^{2k-1}X_{2k-5}+...+q_{3}^{2k-1}X_{3}+q_{1}^{2k-1}X_{1} 
$
\end{center}

\strut \strut

and

\begin{center}
$%
G^{2k}=X_{1}^{2k}=X_{2k}+p_{2k-2}^{2k}X_{2k-2}+p_{2k-4}^{2k}X_{2k-4}+...+p_{2}^{2k}X_{2}+p_{0}^{2k}e, 
$
\end{center}

\strut

where $k\geq 2.$ Also, we have the following recurrence relation, by the
straightforward computation  ;

\strut \strut

\begin{proposition}
Let's fix $k\in \mathbb{N}\,\setminus \,\{1\}.$ Let $q_{i}^{2k-1}$ and $%
p_{j}^{2k}$ ($i=1,3,5,...,2k-1,....$ and $j=0,2,4,...,2k,...$) be given as
before. If $p_{0}^{2}=2N$ and $q_{1}^{3}=(2N-1)+(2N)^{2},$ then we have the
following recurrence relations ;

\strut 

(1) Let

\strut 

$\ \ \ \ \ \ \ \ \
G^{2k-1}=X_{2k-1}+q_{2k-3}^{2k-1}X_{2k-3}+...+q_{3}^{2k-1}X_{3}+q_{1}^{2k-1}X_{1}.
$

Then

\strut 

$\ G^{2k}=X_{2k}+\left( (2N-1)+q_{2k-3}^{2k-1}\right) X_{2k-2}+\left(
(2N-1)q_{2k-3}^{2k-1}+q_{2k-5}^{2k-1}\right) X_{2k-4}$

\strut 

$\ \ \ \ \ \ \ \ \ \ \ \ \ \ \ \ \ +\left(
(2N-1)q_{2k-5}^{2k-1}+q_{2k-7}^{2k-1}\right) X_{2k-6}+$

\strut 

$\ \ \ \ \ \ \ \ \ \ \ \ \ \ \ \ \ +...+\left(
(2N-1)q_{3}^{2k-1}+q_{1}^{2k-1}\right) X_{2}+(2N)q_{1}^{2k-1}e.$

\strut i.e,

\strut 

$\ \ \ \ \ \ \ p_{2k-2}^{2k}=(2N-1)+q_{2k-3}^{2k-1},$ $\ $

$\ \ \ \ \ \ \ p_{2k-4}^{2k}=(2N-1)q_{2k-3}^{2k-1}+q_{2k-5}^{2k-1},$

.....$...,$

\ \ \ \ \ \ \ $p_{2}^{2k}=(2N-1)q_{3}^{2k-1}+q_{1}^{2k-1}$

and

$\ \ \ \ \ \ \ \ p_{0}^{2k}=(2N)q_{1}^{2k-1}.$

\strut 

(2) Let

\strut 

$\ \ \ \ \ \ \ \ \ \ \ \
G^{2k}=X_{2k}+p_{2k-2}^{2k}X_{2k-2}+...+p_{2}^{2k}X_{2}+p_{0}^{2k}e.$

Then

\strut 

$\ G^{2k+1}=X_{2k+1}+\left( (2N-1)+p_{2k-2}^{2k}\right) X_{2k-1}+\left(
(2N-1)p_{2k-2}^{2k}+p_{2k-4}^{2k}\right) X_{2k-3}$

\strut 

$\ \ \ \ \ \ \ \ \ \ \ \ \ \ \ \ \ \ \ \ \ +\left(
(2N-1)p_{2k-4}^{2k}+p_{2k-6}^{2k}\right) X_{2k-5}+$

\strut 

$\ \ \ \ \ \ \ \ \ \ \ \ \ \ \ \ \ \ \ \ \ +...+\left(
(2N-1)p_{4}^{2k}+p_{2}^{2k}\right) X_{3}+\left(
(2N-1)p_{2}^{2k}+p_{0}^{2k}\right) X_{1}.$

i.e,

\strut 

$\ \ \ \ \ \ \ q_{2k-1}^{2k+1}=(2N-1)+p_{2k-2}^{2k},$ \ 

$\ \ \ \ \ \ \ q_{2k-3}^{2k+1}=(2N-1)p_{2k-2}^{2k}+p_{2k-4}^{2k},$

...$...,$ \ 

$\ \ \ \ \ \ \ q_{3}^{2k+1}=(2N-1)p_{4}^{2k}+p_{2}^{2k}$

and

$\ \ \ \ \ \ \ q_{1}^{2k+1}=(2N-1)p_{2}^{2k}+p_{0}^{2k}.$ \ \ \ 

$\square $
\end{proposition}

\strut

\strut

\begin{example}
Suppose that $N=2.$ and let $p_{0}^{2}=4$ and $q_{1}^{3}=3+p_{0}^{2}=3+4=7.$
Put

\strut 

$\ \ \ \ \ \ \
G^{8}=X_{8}+p_{6}^{8}X_{6}+p_{4}^{8}X_{4}+p_{2}^{8}X_{4}+p_{0}^{8}e.$

\strut 

Then, by the previous proposition, we have that

\strut 

$\ p_{6}^{8}=3+q_{5}^{7},$ \ \ $p_{4}^{8}=3q_{5}^{7}+q_{3}^{7},$ \ \ $%
p_{2}^{8}=3q_{3}^{7}+q_{1}^{7}$ \ \ and \ $p_{0}^{8}=4q_{1}^{7}.$

\strut 

Similarly, by the previous proposition,

\strut 

\ \ \ \ $\ q_{5}^{7}=3+p_{4}^{6},$ \ \ \ $q_{3}^{7}=3p_{4}^{6}+p_{2}^{6}$ \
\ \ \ and \ \ $\ q_{1}^{7}=3p_{2}^{6}+p_{0}^{6},$

\strut 

$\ \ \ \ \ \ p_{4}^{6}=3+q_{3}^{5},$ \ \ \ $p_{2}^{6}=3q_{3}^{5}+q_{1}^{5}$
\ \ \ \ \ and \ \ \ \ $p_{0}^{6}=4q_{1}^{5},$

\strut 

$\ \ \ \ \ \ q_{3}^{5}=3+p_{2}^{4}$ \ \ \ \ \ \ \ and \ \ \ \ \ \ $%
q_{1}^{5}=3p_{2}^{4}+p_{2}^{4},$

\strut 

$\ \ \ \ \ \ p_{2}^{4}=3+q_{1}^{3}$ \ \ \ \ \ \ \ \ \ and \ \ \ \ \ \ \ \ \ $%
p_{0}^{4}=4q_{1}^{3},$

\strut 

and

$\ \ \ \ \ \ \ \ \ \ \ q_{1}^{3}=3+p_{0}^{2}=7.$

\strut 

Therefore, combining all information,

\strut 

$\ \ \ \ \ \ G^{8}=X_{8}+22\,X_{6}+202\,X_{4}+744\,X_{2}+1316\,e.$
\end{example}

\strut

We have the following diagram with arrows which mean that

\strut 

\begin{center}
$\swarrow \swarrow $ \ : \ $(2N-1)+[$former term$]$

$\searrow $ \ \ \ \ : \ \ $(2N-1)\cdot [$former term$]$

$\swarrow $ \ \ \ \ : \ \ $\cdot +[$former term$]$
\end{center}

and

\begin{center}
$\searrow \searrow $ \ : \ \ $(2N)\cdot [$former term$].$
\end{center}

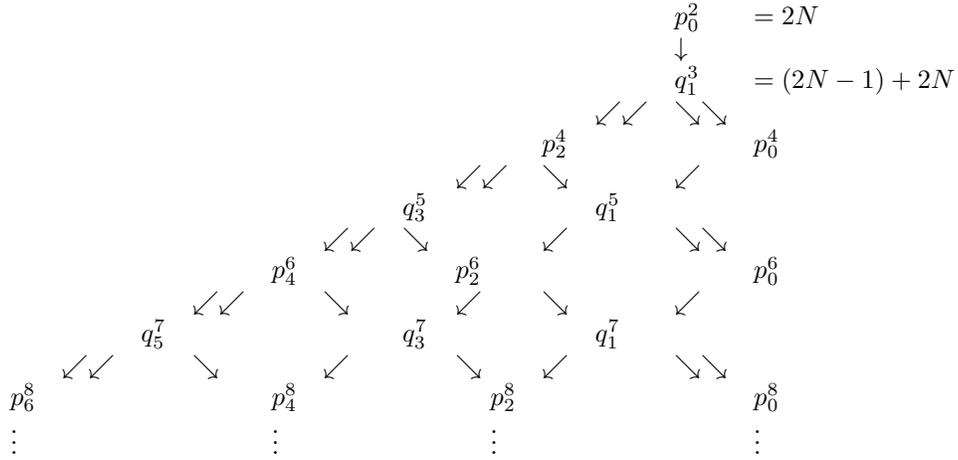
\strut 

\begin{center}
$
\begin{array}{llllllllllll}
&  &  &  &  &  &  &  &  &  & p_{0}^{2} & =2N \\ 
&  &  &  &  &  &  &  &  &  & \downarrow &  \\ 
&  &  &  &  &  &  &  &  &  & q_{1}^{3} & =(2N-1)+2N \\ 
&  &  &  &  &  &  &  &  & \swarrow \swarrow & \searrow \searrow &  \\ 
&  &  &  &  &  &  &  & p_{2}^{4} &  &  & p_{0}^{4} \\ 
&  &  &  &  &  &  & \swarrow \swarrow & \searrow &  & \swarrow &  \\ 
&  &  &  &  &  & q_{3}^{5} &  &  & q_{1}^{5} &  &  \\ 
&  &  &  &  & \swarrow \swarrow & \searrow &  & \swarrow &  & \searrow
\searrow &  \\ 
&  &  &  & p_{4}^{6} &  &  & p_{2}^{6} &  &  &  & p_{0}^{6} \\ 
&  &  & \swarrow \swarrow &  & \searrow &  & \swarrow & \searrow &  & 
\swarrow &  \\ 
&  & q_{5}^{7} &  &  &  & q_{3}^{7} &  &  & q_{1}^{7} &  &  \\ 
& \swarrow \swarrow &  & \searrow &  & \swarrow &  & \searrow & \swarrow & 
& \searrow \searrow &  \\ 
p_{6}^{8} &  &  &  & p_{4}^{8} &  &  & \text{ \ \ \ }p_{2}^{8} &  &  &  & 
p_{0}^{8} \\ 
\vdots &  &  &  & \vdots &  &  & \text{ \ \ \ }\vdots &  &  &  & \vdots
\end{array}
$
\end{center}

\strut 

\begin{quote}
\frame{\textbf{Notation}} From now, we will call the above diagram the%
\textbf{\ recurrence diagram for }$N.$ $\square $
\end{quote}

\strut

For examplet, when $N=2,$ we can compute $p_{4}^{6},$ as follows ;

\strut

\begin{center}
$p_{0}^{2}=4,$ \ \ \ \ \ $\ \ \ q_{1}^{3}=7,$
\end{center}

\strut

\begin{center}
$p_{2}^{4}=3+7=10,$ \ \ \ \ \ $p_{0}^{4}=28,$
\end{center}

\strut

\begin{center}
$q_{3}^{5}=3+10=13,$ \ $\ q_{1}^{5}=3\cdot 10+28=58.$
\end{center}

\strut

and hence $p_{4}^{6}=3+13=16.$

\strut \strut 

By the recurrence diagram for $N$, we can get that ;

\strut \strut 

\begin{theorem}
Let $G\in \left( A,\tau \right) $ be the generating operator. Then the
moment series of $G$ is

\strut 

$\ \ \ \ \ \ \ \ \ \tau \left( G^{n}\right) =\left\{ 
\begin{array}{lll}
0 &  & \text{if }n\text{ is odd} \\ 
&  &  \\ 
p_{0}^{n} &  & \text{if }n\text{ is even,}
\end{array}
\right. $

\strut 

for all $n\in \mathbb{N}.$
\end{theorem}

\strut

\begin{proof}
Assume that $n$ is odd. Then

\strut

$\ \ \ \ \ \
G^{n}=X_{n}+q_{n-2}^{n}X_{n-2}+...+q_{3}^{n}X_{3}+q_{1}^{n}X_{1}.$

\strut

So, $G^{n}$ does not contain the $e$-terms. Therefore,

\strut

$\ \ \ \tau \left( G^{n}\right) =\tau \left(
X_{n}+q_{n-2}^{n}X_{n-2}+...+q_{3}^{n}X_{3}+q_{1}^{n}X_{1}\right) =0.$

\strut

Assume that $n$ is even. Then

\strut

\ $\ \ \ \ G^{n}=X_{n}+p_{n-2}^{n}X_{n-2}+...+p_{2}^{n}X_{2}+p_{0}^{n}e.$

\strut

So, we have that

\strut

$\ \ \ \tau (G^{n})=\tau \left(
X_{n}+p_{n-2}^{n}X_{n-2}+...+p_{2}^{n}X_{2}+p_{0}^{n}e\right) =p_{0}^{n}.$
\end{proof}

\strut

Remark that the $n$-th moments of the generating operator in $(A,\tau )$ is
totally depending on the recurrence diagram for $N$.\strut

\strut 

The moments of the generating operator $G$ of $A=L(F_{N})$ is well-know. But
the above theorem provides the method how to get the moments of generating
operator $G$ by using the recurrence diagram. 

\strut 

Recall that Nica and Speicher defined the even random variable in a $*$%
-probability space. Let $(B,\tau _{0})$ be a $*$-probability space, where $%
\tau _{0}$ $:$ $B\rightarrow $ $\mathbb{C}$ is a linear functional
satisfying that $\tau _{0}\left( b^{*}\right) =\overline{\tau _{0}(b)},$ for
all $b\in B,$ and let $b\in (B,\tau _{0})$ be a random variable. We say that
the random variable $b\in (B,\tau _{0})$ is even if it is self-adjoint and
it satisfies the following moment relation ;

\strut

\begin{center}
$\tau _{0}\left( b^{n}\right) =0,$ whenever $n$ is odd.
\end{center}

\strut 

By the previous theorm, the generating operator $G$ of $A=L(F_{N})$ is an
even element (for all $N\in \Bbb{N}$).

\strut \strut \strut \strut \strut 

\begin{corollary}
Let $G\in \left( A,\tau \right) $ be the generating operator. Then $G$ is
even in $\left( A,\tau \right) $. $\square $
\end{corollary}

\strut 

Also, by the previous theorem, we can recompute the well-known moment series
of the generating operator $G$ of $L(F_{N})$ ;

\strut 

\begin{corollary}
Let $G\in (A,\tau )$ be the generating operator. Then the operator $G$ has
its moment series,

\strut 

$\ \ \ \ \ \ \ \ \ \ \ M_{G}(z)=\sum_{n=1}^{\infty }p_{0}^{2n}\,z^{2n}\in
\Theta _{1},$

\strut 

where $p_{0}^{2n}$ are the given numbers in the recurrence diagram. \ $%
\square $
\end{corollary}

\strut \ \strut

\strut

\strut

\section{Operator-valued Moment Series of the Generating Operator of $%
L(F_{2})$}

\strut

\strut

\strut

Let $M_{0}\subset M$ be von Neumann algebras with $1_{M_{0}}=1_{M}$ and let $%
\varphi :M\rightarrow M_{0}$ be the conditional expectation. i.e, the
surjective bimodule map $\varphi $ satisfies that

\strut

\begin{center}
$\varphi (m_{0})=m_{0},$ for all $m_{0}\in M,$
\end{center}

\strut

\begin{center}
$\varphi (m_{0}mm_{0}^{\prime })=m_{0}\cdot \varphi (m)\cdot m_{0}^{\prime
}, $
\end{center}

\strut

for all $m_{0},m_{0}^{\prime }\in M_{0},$ $m\in M,$ and

\strut

\begin{center}
$\varphi (m^{*})=\varphi (m)^{*},$ for all $m\in M.$
\end{center}

\strut

Then the algebraic pair $(M,\varphi )$ is a $W^{*}$-probability space with
amalgamation over $M_{0}$ (See [1]). If $m\in (M,\varphi ),$ then we will
call $m$ an operator-valued or ($M_{0}$-valued) random variable.

\strut

\begin{definition}
Let $(M,\varphi )$ be a $W^{*}$-probability space over $M_{0}$ and let $m\in
(M,\varphi )$ be a $M_{0}$-valued random variable. Define the $n$-th
operator-valued moment of $m$ by

\strut 

$\ \ \ \ \ \ \ \ \ E\left( (m_{1}m)(m_{2}m)...(m_{n}m)\right) ,$

\strut 

for all $n\in \mathbb{N},$ where $m_{1},...,m_{n}\in M_{0}$ are arbitrary.
When $m_{1}=...=m_{n}=1_{M_{0}},$ for all $n\in \mathbb{N},$ we say that the 
$n$-th operator-valued moments of $m$ are trivial. i.e, the $n$-th \textbf{%
trivial} $M_{0}$-valued moments of $m\in (M,\varphi )$ are $E(m^{n}),$ for
all $n\in \mathbb{N}.$ We will say that the operator-valued formal series

\strut 

$\ \ \ \ \ \ \ \ \ M_{m}^{t}(z)=\sum_{n=1}^{\infty }E(m^{n})\,z^{n}\in
M_{0}[[z]]$

\strut 

is the t\textbf{rivial }$M_{0}$\textbf{-valued moment series of }$m\in
(M,\varphi ),$ where $z$ is the indeterminent.
\end{definition}

\strut \strut

In this chapter, by using the recurrence diagram for $N=2,$ we will compute
the trivial operator-valued moment series of the given generating operator

\strut

\begin{center}
$G=a+b+a^{-1}+b^{-1}$
\end{center}

\strut

over the commutator-group von Neumann algebra $L(K),$ where $F_{2}=\,<a,b>$
is a free group and $K=\,<h=aba^{-1}b^{-1}>$ is the commutator group of $a$
and $b$.

\strut

First, we will define the conditional expectation $E:L(F_{2})\rightarrow L(K)
$ by

\strut

(2.1) $\ \ \ \ \ \ \ \ \ \ \ \ E\left( \underset{g\in F_{2}}{\sum }\alpha
_{g}g\right) =\underset{k\in K}{\sum }\alpha _{k}k,$

\strut

for all $\underset{g\in F_{2}}{\sum }\alpha _{g}g\in \left(
L(F_{2}),E\right) .$

\strut

Then we can construct the $W^{*}$-probability space $\left(
L(F_{2}),E\right) $ over its $W^{*}$-subalgebra $L(K).$ Notice that to find
the conditional expectational value of the operator-valued random variable $%
x\in \left( L(F_{2}),E\right) $ is to find the $h^{k}$-terms of the
operator-valued random variables $x$, for $k\in \mathbb{Z},$ where

\strut 

\begin{center}
$h=aba^{-1}b^{-1}$ \ \ and \  $h^{-1}=bab^{-1}a^{-1}.$
\end{center}

\strut

First, let us provide the recurrence diagram for $N=2$ ;

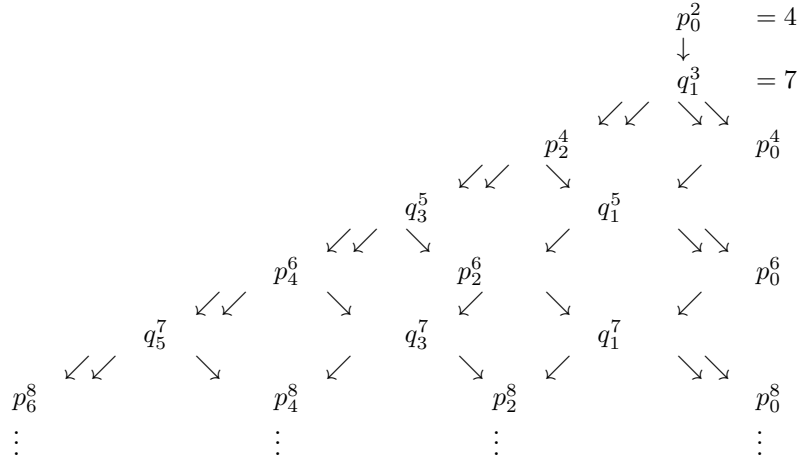
\strut

\strut

\begin{center}
$
\begin{array}{llllllllllll}
&  &  &  &  &  &  &  &  &  & p_{0}^{2} & =4 \\ 
&  &  &  &  &  &  &  &  &  & \downarrow &  \\ 
&  &  &  &  &  &  &  &  &  & q_{1}^{3} & =7 \\ 
&  &  &  &  &  &  &  &  & \swarrow \swarrow & \searrow \searrow &  \\ 
&  &  &  &  &  &  &  & p_{2}^{4} &  &  & p_{0}^{4} \\ 
&  &  &  &  &  &  & \swarrow \swarrow & \searrow &  & \swarrow &  \\ 
&  &  &  &  &  & q_{3}^{5} &  &  & q_{1}^{5} &  &  \\ 
&  &  &  &  & \swarrow \swarrow & \searrow &  & \swarrow &  & \searrow
\searrow &  \\ 
&  &  &  & p_{4}^{6} &  &  & p_{2}^{6} &  &  &  & p_{0}^{6} \\ 
&  &  & \swarrow \swarrow &  & \searrow &  & \swarrow & \searrow &  & 
\swarrow &  \\ 
&  & q_{5}^{7} &  &  &  & q_{3}^{7} &  &  & q_{1}^{7} &  &  \\ 
& \swarrow \swarrow &  & \searrow &  & \swarrow &  & \searrow & \swarrow & 
& \searrow \searrow &  \\ 
p_{6}^{8} &  &  &  & p_{4}^{8} &  &  & \text{ \ \ \ }p_{2}^{8} &  &  &  & 
p_{0}^{8} \\ 
\vdots &  &  &  & \vdots &  &  & \text{ \ \ \ }\vdots &  &  &  & \vdots
\end{array}
$
\end{center}

\strut

where

\begin{center}
$\swarrow \swarrow $ \ : \ $3+[$former term$]$

$\searrow $ \ \ \ \ : \ \ $3\cdot [$former term$]$

$\swarrow $ \ \ \ \ : \ \ $\cdot +[$former term$]$
\end{center}

and

\begin{center}
$\searrow \searrow $ \ : \ \ $4\cdot [$former term$].$
\end{center}

\strut

By the above recurrence diagram for $N=2,$ we have that if

$\strut $

\begin{center}
$p_{0}^{2}=4$ \ \ \ \ and $\ \ \ q_{1}^{3}=3+p_{0}^{2}=7,$
\end{center}

\strut

and if we put $X_{n}=\underset{\left| w\right| =n}{\sum }w,$ as the sum of
all words with length $n$ in $\{a,b,a^{-1},b^{-1}\},$ for $n\in \mathbb{N},$
then we have the following recurrence relations (1) and (2) ;

\strut

(1) If

\strut \strut

\begin{center}
$%
G^{2k-1}=X_{2k-1}+q_{2k-3}^{2k-1}X_{2k-3}+...+q_{3}^{2k-1}X_{3}+q_{1}^{2k-1}X_{1}. 
$
\end{center}

then

\strut

$\ \ \ \ G^{2k}=X_{2k}+\left( 3+q_{2k-3}^{2k-1}\right) X_{2k-2}+\left(
3q_{2k-3}^{2k-1}+q_{2k-5}^{2k-1}\right) X_{2k-4}$

\strut

$\ \ \ \ \ \ \ \ \ \ \ \ \ \ \ +\left(
3q_{2k-5}^{2k-1}+q_{2k-7}^{2k-1}\right) X_{2k-6}$

\strut

$\ \ \ \ \ \ \ \ \ \ \ \ \ \ \ +...+\left( 3q_{3}^{2k-1}+q_{1}^{2k-1}\right)
X_{2}+4q_{1}^{2k-1}e.$

\strut where

\strut

\begin{center}
$p_{2k-2}^{2k}=3+q_{2k-3}^{2k-1},$ $\ \ \ \ \
p_{2k-4}^{2k}=3q_{2k-3}^{2k-1}+q_{2k-5}^{2k-1},$
\end{center}

$\strut $

\begin{center}
$...,$ $p_{2}^{2k}=3q_{3}^{2k-1}+q_{1}^{2k-1}$ \ \ \ and \ \ \ $%
p_{0}^{2k}=4q_{1}^{2k-1},$
\end{center}

\strut

by the recurrence diagram.

\strut

(2) If

\begin{center}
$G^{2k}=X_{2k}+p_{2k-2}^{2k}X_{2k-2}+...+p_{2}^{2k}X_{2}+p_{0}^{2k}e.$
\end{center}

Then

\strut

$\ \ \ \ \ \ G^{2k+1}=X_{2k+1}+\left( 3+p_{2k-2}^{2k}\right) X_{2k-1}+\left(
3p_{2k-2}^{2k}+p_{2k-4}^{2k}\right) X_{2k-3}$

\strut

$\ \ \ \ \ \ \ \ \ \ \ \ \ \ \ \ \ \ \ \ \ +\left(
3p_{2k-4}^{2k}+p_{2k-6}^{2k}\right) X_{2k-5}+$

\strut

$\ \ \ \ \ \ \ \ \ \ \ \ \ \ \ \ \ \ \ \ \ +...+\left(
3p_{4}^{2k}+p_{2}^{2k}\right) X_{3}+\left( 3p_{2}^{2k}+p_{0}^{2k}\right)
X_{1}.$

where

\strut

\begin{center}
$q_{2k-1}^{2k+1}=3+p_{2k-2}^{2k},$ \ \ \ \ \ $%
q_{2k-3}^{2k+1}=3p_{2k-2}^{2k}+p_{2k-4}^{2k},$
\end{center}

$\strut $

\begin{center}
$...,$ \ $q_{3}^{2k+1}=3p_{4}^{2k}+p_{2}^{2k}$ \ \ \ and \ \ $%
q_{1}^{2k+1}=3p_{2}^{2k}+p_{0}^{2k},$
\end{center}

\strut

by the recurrence diagram.

\strut \strut

Note that $h$ and $h^{-1}$ are words with their length 4. Therefore, $X_{4k}$
contains $h^{k}$-terms and $h^{-k}$-terms, for all $k\in \mathbb{N\cup \{}0%
\mathbb{\}}$ ! Thus we can compute the trivial operator-valued moments of
the operator $G$ as follows ;

\strut

\begin{theorem}
Fix $k\in \mathbb{N}$ and Let $G\in \left( L(F_{2}),E\right) $ be the
generating operator of $L(F_{2})$. Then

\strut 

(1) $\ E(G^{k})=0_{L(K)},$ if $k$ is odd.

\strut 

(2) $\ E\left( G^{4k}\right) =\left( h^{k}+h^{-k}\right)
+\sum_{j=1}^{k-1}p_{4k-4j}^{4k}\left( h^{k-j}+h^{-(k-j)}\right)
+p_{0}^{4k}h^{0},$

\strut 

where $p_{0}^{4}=28.$

\strut 

(3) If $4\nmid 2k,$ in the sense that $2k$ is not a multiple by 4, then

\strut 

$\ \ \ \ \ \ E(G^{2k})=\sum_{j=1}^{k-1}p_{(2k-2)-4j}^{2k}\left( h^{\frac{k-1%
}{2}-2j}+h^{-(\frac{k-1}{2}-2j)}\right) +p_{0}^{2k}h^{0},$

\strut 

where $p_{0}^{2}=4.$
\end{theorem}

\strut

\begin{proof}
(1) Suppose that $k$ is odd. Then $G^{k}$ does not have the words with
length $4p,$ for some $p\in \mathbb{N},$ by the recurrence diagram for $N=2,$
since $G^{k}$ does not have the $X_{4n}$-terms, for $n\in \mathbb{N},$ $%
4n<k. $ This shows that there's no $h^{n}$-terms and $h^{-n}$-terms in $%
G^{k},$ where $n$ is previousely given such that $2n<k.$ Therefore, all odd
trivial operator-valued moments of $G$ vanish.

\strut

(2) By the straightforward computation using the recurrence diagram, we have
that

\strut

$\ E\left( G^{4k}\right) $

\strut

$\ \ =E\left(
X_{4k}+p_{4k-2}^{4k}X_{4k-2}+p_{4k-4}^{4k}X_{4k-4}+...+p_{4}^{4k}X_{4}+p_{2}^{4k}X_{2}+p_{0}^{4k}h^{0}\right) 
$

\strut

\ \ (2.2)\strut

\strut

$\ \ =E(X_{4k})+p_{4k-2}^{4k}E(X_{4k-2})+p_{4k-4}^{4k}E(X_{4k-4})+$

\strut

$\ \ \ \ \ \ \ \ \ \ \ \ \ \ \
...+p_{4}^{4k}E(X_{4})+p_{2}^{4k}E(X_{2})+p_{0}^{4k}h^{0}.$

\strut

Since $h^{p}$ and $h^{-p}$ terms are in $X_{4p},$ for any $p\in \mathbb{N}%
\cup \{0\},$ the formular (2.2) is

\strut

\ \ (2.3)

\strut

$\ \ \ \ \
E(X_{4k})+p_{4k-4}^{4k}E(X_{4k-4})+...+p_{4}^{4k}E(X_{4})+p_{0}^{4k}h^{0}$

\strut

$\ \ \ \ =\left( h^{k}+h^{-k}\right)
+p_{4k-4}^{4k}(h^{k-1}+h^{-(k-1)})+...+p_{4}^{4k}(h+h^{-1})+p_{0}^{4k}h^{0}.$

\strut

(3) If $4\nmid 2k,$ then $k=1,3,5,....$. If $k=1,$ then we have that ;

\strut

$\ \ \ \ \ \ \ \ \ \ \ \ \ \ \ \ \ E(G^{2})=E\left( X_{2}+4h^{0}\right)
=4h^{0}.$

\strut

If $k\neq 1$ is odd, then

\strut

$\ E(G^{2k})$

\strut

$\ \
=E(X_{2k}+p_{2k-2}^{2k}X_{2k-2}+p_{2k-4}^{2k}X_{2k-4}+p_{2k-6}^{2k}X_{2k-6}+$

\strut

$\ \ \ \ \ \ \ \ \ \ \ \ \ \ \ \ \ \ \ \ \ \ \ \ \ \ \ \
...+p_{4}^{2k}X_{4}+p_{2}^{2k}X_{2}+p_{0}^{2k}h^{0})$

\strut

$\ \
=E(X_{2k})+p_{2k-2}^{2k}E(X_{2k-2})+p_{2k-4}^{2k}E(X_{2k-4})+p_{2k-6}^{2k}E(X_{2k-6})+ 
$

\strut

$\ \ \ \ \ \ \ \ \ \ \ \ \ \ \ \ \ \ \ \ \ \ \ \ \ \ \ \
...+p_{4}^{2k}E(X_{4})+p_{2}^{2k}E(X_{2})+p_{0}^{2k}h^{0}$

\strut

$\ =0_{B}+p_{2k-2}^{2k}\left( h^{k-1}+h^{-(k-1)}\right)
+0_{B}+p_{2k-6}^{2k}\left( h^{k-3}+h^{-(k-3)}\right) +$

\strut

$\ \ \ \ \ \ \ \ \ \ \ \ \ \ \ \ \ \ \ \ \ \ \ \ \ \ \ \ \ \
...+p_{4}^{2k}(h+h^{-1})+0_{B}+p_{0}^{2k}h^{0},$

\strut

since $X_{2k-2},$ $X_{2k-6},...,X_{4}$ contain $h^{p}$-terms and $h^{-p}$%
-terms, for $p\in \mathbb{N}\cup \{0\}.$
\end{proof}

\strut

By the previous trivial operator-valued moments of the generating operator $%
G $ of $L(F_{2}),$ we have the following result ;

\strut

\begin{corollary}
Let $\left( L(F_{2}),E\right) $ be the $W^{*}$-probability spcae over $L(K)$
and let $G\in \left( L(F_{2}),E\right) $ be the generating operator of $%
L(F_{2})$. Then the trivial operator-valued moment series of $G$ is

\strut 

$\ \ \ \ \ \ \ \ \ \ \ M_{G}(z)=\sum_{n=1}^{\infty }b_{2n}\,z^{n}\in
L(K)[[z]],$

\strut 

where

\strut 

$\ \ \ \ \ \ b_{4n}=\left( h^{n}+h^{-n}\right)
+\sum_{j=1}^{n-1}p_{4n-4j}^{4n}\left( h^{n-j}+h^{-(n-j)}\right)
+p_{0}^{4n}h^{0}$

and

$\ \ \ \ \ \ \ \ b_{2k}=\sum_{j=1}^{k-1}p_{(2k-2)-4j}^{2k}\left( h^{\frac{k-1%
}{2}-2j}+h^{-(\frac{k-1}{2}-2j)}\right) +p_{0}^{2k}h^{0},$

\strut 

where $4\nmid 2k,$ for all $n,k\in \mathbb{N}.$ $\square $
\end{corollary}

\strut

\begin{remark}
Suppose we have the free group factor $L(F_{N}),$ where $N\in \mathbb{N}.$
Then we can extend the above result for the general $N.$ i.e, we can take $%
F_{N}=~<g_{1},...,g_{N}>$ \ and \ $K=~<h>,$ where

\strut 

$\ \ \ \ \ \ \ \ \ k=g_{1}\cdot \cdot \cdot g_{N}\cdot g_{1}^{-1}\cdot \cdot
\cdot g_{N}^{-1}.$

\strut 

By defining the canonical conditional expectation $E:L(F_{N})\rightarrow
L(K),$ we can construct the $W^{*}$-probability space with amalgamation over 
$L(K)$. Similar to the case when $N=2,$ by using the recurrence diagram for $%
N,$ we can compute the trivial operator-valued moments of the generating
operator

\strut 

$\ \ \ \ \ \ \ G_{N}=g_{1}+...+g_{N}+g_{1}^{-1}+...+g_{N}^{-1}.$

\strut 

When $N=2,$ to find the trivial operator-valued moments of $G_{2}$ is to
find the $h^{k}$-terms of $G_{2}^{n},$ for $k\in \mathbb{Z},$ $n\in \mathbb{N%
}.$ But in the general case, to find the trivial operator-valued moments of $%
G_{N}$ is to find the $k^{p}$-terms of $G_{N}^{q},$ for $p\in \mathbb{Z},$ $%
q\in \mathbb{N}.$ So, we have to choose the $X_{(2N)p}$-terms, for all $p\in %
\mathbb{N}$, containing $k^{\pm p}$-terms !
\end{remark}

\strut

\strut

\strut

\strut \textbf{References}

\bigskip

\strut

{\small [1] \ R. Speicher, Combinatorial Theory of the Free Product with
Amalgamation and Operator-Valued Free Probability Theory, AMS Mem, Vol 132 ,
Num 627 , (1998).}

{\small [2] \ \ A. Nica, R-transform in Free Probability, IHP course note.}

{\small [3] \ \ R. Speicher, Combinatorics of Free Probability Theory IHP
course note.}

{\small [4] \ \ A. Nica, D. Shlyakhtenko and R. Speicher, R-cyclic Families
of Matrices in Free Probability, J. of Funct Anal, 188 (2002),
227-271.\strut }

{\small [5] \ \ A. Nica and R. Speicher, R-diagonal Pair-A Common Approach
to Haar Unitaries and Circular Elements, (1995), Preprint.}

{\small [6] \ \ D. Shlyakhtenko, Some Applications of Freeness with
Amalgamation, J. Reine Angew. Math, 500 (1998), 191-212.\strut }

{\small [7] \ \ A. Nica, D. Shlyakhtenko and R. Speicher, R-diagonal
Elements and Freeness with Amalgamation, Canad. J. Math. Vol 53, Num 2,
(2001) 355-381.\strut }

{\small [8] \ \ A. Nica, R-transforms of Free Joint Distributions and
Non-crossing Partitions, J. of Func. Anal, 135 (1996), 271-296.\strut }

{\small [9] \ \ D.Voiculescu, K. Dykemma and A. Nica, Free Random Variables,
CRM Monograph Series Vol 1 (1992).\strut }

{\small [10] D. Voiculescu, Operations on Certain Non-commuting
Operator-Valued Random Variables, Ast\'{e}risque, 232 (1995), 243-275.\strut 
}

{\small [11] D. Shlyakhtenko, A-Valued Semicircular Systems, J. of Funct
Anal, 166 (1999), 1-47.\strut }

{\small [12] I. Cho, The Moment Series and The R-transform of the Generating
Operator of }$L(F_{N}),${\small \ (2003), Preprint.\strut }

{\small [13] I. Cho, The Moment Series of the Generating Operator of }$%
L(F_{2})*_{L(K)}L(F_{2})${\small , (2003), Preprint}

{\small [14] I. Cho, An Example of Moment Series under the Compatibility,
(2003), Preprint}

{\small [15] F. Radulescu, Singularity of the Radial Subalgebra of }$%
L(F_{N}) ${\small \ and the Puk\'{a}nszky Invariant, Pacific J. of Math,
vol. 151, No 2 (1991)\strut , 297-306.\strut \strut }

\strut

\end{document}